\documentclass[11pt]{article}
\usepackage{amssymb,amsfonts}
\begin{document}

\begin{center}
\Large On movable singularities of Garnier systems
\end{center}

\centerline{R.\,R.\,Gontsov}

\begin{abstract}
We study movable singularities of Garnier systems using the
connection of the latter with Schlesinger isomonodromic
deformations of Fuchsian systems.
\end{abstract}
\medskip

{\bf \S1. What is Painlev\'e VI equations and Garnier systems?}
\medskip

We start with the Painlev\'e VI ($\rm P_{VI}$) equation
\begin{eqnarray}\label{PVI}
\frac{d^2u}{dt^2}&=&\frac12\left(\frac1u+\frac1{u-1}+\frac1{u-t}\right)\left(\frac{du}{dt}\right)^2-
\left(\frac1t+\frac1{t-1}+\frac1{u-t}\right)\frac{du}{dt}+\nonumber \\
& &+\frac{u(u-1)(u-t)}{t^2(t-1)^2}\left(\alpha+\beta\,\frac
t{u^2}+\gamma\,\frac{t-1}{(u-1)^2}+\delta\,\frac{t(t-1)}{(u-t)^2}\right),
\end{eqnarray}
the second order ODE for a complex function $u(t)$, where
$\alpha,\beta,\gamma,\delta\in{\mathbb C}$ are constants.

However, simply giving the explicit equation seems to be the least
helpful introduction to it. For the purposes of this paper it is
more convenient to look at $\rm P_{VI}$ as at
\begin{itemize}
\item the equation for an apparent (fifth) singularity of
isomonodromic family of second order scalar Fuchsian equations
with the four singularities $t,0,1,\infty$;

\item the most general second order ODE with the Painlev\'e
property;

\item the equation controlling isomonodromic deformations of
certain rank 2 Fuchsian systems with the four singularities
$t,0,1,\infty$.
\end{itemize}

Let us recall the first two viewpoints in more details (the last
one will appear in \S3).

The monodromy of a linear differential equation
\begin{eqnarray}\label{urav}
\frac{d^pu}{dz^p}+b_1(z)\frac{d^{p-1}u}{dz^{p-1}}+\ldots+b_p(z)u=0
\end{eqnarray}
with singularities $a_1,\ldots,a_n\in\overline{\mathbb C}$ (which
are the poles of the coefficients) can be defined as follows. In a
neighbourhood of a non-singular point $z_0$ we consider a basis
$(u_1,\ldots,u_p)$ in the solution space of the equation
(\ref{urav}). Analytic continuations of the functions
$u_1(z),\ldots,u_p(z)$ along an arbitrary loop $\gamma$ outgoing
from $z_0$ and lying in $\overline{\mathbb C}\setminus
\{a_1,\ldots,a_n\}$ transform the basis $(u_1,\ldots,u_p)$ into a
(in general case different) basis $(\tilde u_1,\ldots,\tilde
u_p)$. The two bases are related by means of a non-singular
transition matrix $G_{\gamma}$ corresponding to the loop $\gamma$:
$$
(u_1,\ldots,u_p)=(\tilde u_1,\ldots,\tilde u_p)G_{\gamma}.
$$
The map $[\gamma]\mapsto G_{\gamma}$ (which depends only on the
homotopy class $[\gamma]$ of the loop $\gamma$) defines the
representation
$$
\chi: \pi_1(\overline{\mathbb C}\setminus\{a_1,\ldots,a_n\},z_0)
\longrightarrow{\rm GL}(p,\mathbb C)
$$
of the fundamental group of the space $\overline{\mathbb
C}\setminus \{a_1,\ldots, a_n\}$ in the space of non-singular
complex matrices of size $p$. This representation is called the
{\it monodromy} of the equation (\ref{urav}).

A singular point $a_i$ of the equation (\ref{urav}) is said to be
{\it regular} if any solution of the equation has a polynomial
(with respect to $1/|z-a_i|$) growth near $a_i$. Linear
differential equations with regular singular points only are
called {\it Fuchsian}.

A.\,Poincar\'e \cite{Po} has established that the number of
parameters determining a Fuchsian equation of order $p$ with $n$
singular points is less than the dimension of the space of
representations $\chi$, if $p>2, n>2$ or $p=2, n>3$ (see also
\cite{AB}, pp. 158--159). Hence in the construction of a Fuchsian
equation with the given singularities and monodromy there arise
so-called {\it apparent} singularities, at which the coefficients
of the equation have poles but the solutions are single-valued
meromorphic functions. In the case $p=2$, $n=4$
($\{a_1,a_2,a_3,a_4\}=\{t,0,1,\infty\}$) the number of such
singularities equals one. If we move a little the singularity
$z=t$ so that the monodromy of the equation preserves (this is an
{\it isomonodromy} property which is defined precisely in the next
paragraph), the apparent (fifth) singularity $w(t)$ will move
satisfying $\rm P_{VI}$ (this was first obtained by R.\,Fuchs
\cite{Fu1}).

The equation (\ref{PVI}) has three fixed singular points -- $0, 1,
\infty$. Its movable singularities (which depend on the initial
conditions) can be poles only. In other words, any local solution
of the equation defined in a neighbourhood of $t_0\ne0,1,\infty$
can be extended to a meromorphic function on the universal cover
of $\overline{\mathbb C}\setminus\{0,1,\infty\}$. This is the {\it
Painlev\'e property}. The statement on movable poles of the
equation (\ref{PVI}) is the following. In the case $\alpha\ne0$
they can be simple only, and in the case $\alpha=0$ their orders
do not exceed two (see, for instance, \cite{GL}, Ch. VI, \S6).

Extending the first of the above three viewpoints to general case
of $n+3$ singularities $a_1,\ldots,a_n,0,1,\infty$, R.\,Garnier
\cite{Gar} has obtained the system ${\cal G}_n(\theta)$ depending
on $n+3$ complex parameters $\theta_1,\ldots,\theta_{n+2},
\theta_{\infty}$. This is a completely integrable system of
non-linear partial differential equations of second order. Later
it was written down by K.\,Okamoto \cite{Ok} in an equivalent
Hamiltonian form
\begin{eqnarray}\label{garnier}
\frac{\partial u_i}{\partial a_j}=\frac{\partial H_j}{\partial
v_i}, \qquad \frac{\partial v_i}{\partial a_j}= -\frac{\partial
H_j}{\partial u_i}, \qquad i,j=1,\ldots,n,
\end{eqnarray}
with certain Hamiltonians $H_i=H_i(a,u,v,\theta)$ rationally
depending on $a=(a_1,\ldots,a_n)$, $u=(u_1,\ldots,u_n)$,
$v=(v_1,\ldots,v_n)$, $\theta=(\theta_1,\ldots,\theta_{n+2},
\theta_{\infty})$. (Here $u_1(a),\ldots,u_n(a)$ are apparent
singular points of a certain isomonodromic family of second order
Fuchsian equations with singularities
$a_1,\ldots,a_n,0,1,\infty$.) In the case $n=1$ the Garnier system
${\cal G}_1(\theta_1,\theta_2, \theta_3,\theta_{\infty})$ is an
equivalent (Hamiltonian) form of ${\rm P_{VI}}$ (\ref{PVI}), where
$$
\alpha=\frac12\theta_{\infty}^2, \quad \beta=-\frac12\theta_2^2,
\quad \gamma=\frac12\theta_3^2, \quad
\delta=\frac12(1-\theta_1^2).
$$

For $n>1$ the Garnier system generically does not satisfy the
Painlev\'e property. However, due to Garnier's theorem, the
elementary symmetric polynomials $\sigma_i(u_1(a),\ldots,u_n(a))$,
depending on local solutions of the Garnier system, extend to
meromorphic functions $F_i(a)$ on the universal cover $Z'$ of the
space $({\mathbb C}\setminus\{0,1\})^n\setminus\bigcup_{i\ne j}
\{a_i=a_j\}$. Our addition to this theorem consists in some
estimates for orders of irreducible components of the polar loci
of the functions $F_i$ (Theorem 2, Proposition 1).
\\
\\

{\bf \S2. Isomonodromic deformations of Fuchsian systems}
\medskip

Let us include a Fuchsian system
\begin{eqnarray}\label{syst0}
\frac{dy}{dz}=\left(\sum_{i=1}^n\frac{B^0_i}{z-a^0_i}\right)y,
\qquad B_i^0\in{\rm Mat}(p,{\mathbb C}), \qquad
\sum_{i=1}^nB^0_i=0,
\end{eqnarray}
of $p$ equations with singularities $a_1^0,\ldots,a_n^0$ into a
family
\begin{eqnarray}\label{fam}
\frac{dy}{dz}=\left(\sum_{i=1}^n\frac{B_i(a)}{z-a_i}\right)y,
\qquad B_i(a^0)=B^0_i, \qquad \sum_{i=1}^nB_i(a)=0,
\end{eqnarray}
of Fuchsian systems holomorphically depending on the parameter
$a=(a_1,\ldots,a_n)\in D(a^0)$, where $D(a^0)$ is a disk of small
radius centered at the point $a^0=(a_1^0,\ldots,a_n^0)$ of the
space ${\mathbb C}^n\setminus\bigcup_{i\ne j}\{a_i=a_j\}$.

One says that the family (\ref{fam}) is {\it isomonodromic} (or it
is an {\it isomonodromic deformation} of the system
(\ref{syst0})), if for all $a\in D(a^0)$ the monodromies
$$
\chi: \pi_1(\overline{\mathbb C}\setminus\{a_1,\ldots,a_n\})
\longrightarrow{\rm GL}(p,{\mathbb C})
$$
of the corresponding systems are the same\footnote{Under small
variations of the parameter $a$ there exist canonical isomorphisms
of the fundamental groups $\pi_1(\overline{\mathbb C}
\setminus\{a_1,\ldots,a_n\})$ and $\pi_1(\overline{\mathbb C}
\setminus\{a_1^0,\ldots,a_n^0\})$ generating canonical
isomorphisms
$$
{\rm Hom}\left(\pi_1(\overline{\mathbb C}
\setminus\{a_1,\ldots,a_n\}), {\rm GL}(p,{\mathbb C})\right)/{\rm
GL}(p,{\mathbb C})\cong{\rm Hom}\left(\pi_1(\overline{\mathbb
C}\setminus\{a_1^0,\ldots,a_n^0\}),{\rm GL}(p,{\mathbb C})\right)/
{\rm GL}(p,{\mathbb C})
$$
of the spaces of conjugacy classes of representations for the
above fundamental groups; this allows one to compare $\chi$ for
various $a\in D(a^0)$.}. This means that for every value $a$ there
exists a fundamental matrix $Y(z, a)$ of the corresponding system
(\ref{fam}) that has the same monodromy for all $a\in D(a^0)$.
This matrix $Y(z, a)$ is called an {\it isomonodromic fundamental
matrix}.

Is it always possible to include the system (\ref{syst0}) into an
isomonodromic family of Fuchsian systems? The answer is positive.
For instance, if the matrices $B_i(a)$ satisfy the {\it
Schlesinger equation} \cite{Sch}
$$
dB_i(a)=-\sum_{j=1, j\ne i}^n\frac{[B_i(a), B_j(a)]}
{a_i-a_j}\,d(a_i-a_j),
$$
then the family (\ref{fam}) is isomonodromic (in this case it is
called the {\it Schlesinger isomonodromic family}).

Due to Malgrange's theorem \cite{Ma}, for arbitrary initial
conditions $B_i(a^0)=B_i^0$ the Schlesinger equation has the
unique solution $\{B_1(a),\ldots,B_n(a)\}$ in some disk $D(a^0)$,
and the matrices $B_i(a)$ can be extended to the universal cover
$Z$ of the space ${\mathbb C}^n\setminus \bigcup_{i\ne
j}\{a_i=a_j\}$ as meromorphic functions. Thus, the Schlesinger
equation satisfies the Painlev\'e property.

Recall, that a function $f$ is {\it meromorphic} on $Z$, if it is
holomorphic on $Z\setminus P$, can not be extended to $P$
holomorphically and is presented as a quotient
$f(a)=\varphi(a)/\psi(a)$ of holomorphic functions in a
neighbourhood of every point $a^0\in P$ (hence, $\psi(a^0)=0$).
Thus, $P\subset Z$ is an analytic set of codimension one (it is
defined locally by the equation $\psi(a)=0$), which is called the
{\it polar locus} of the meromorphic function $f$. The points of
this set is divided into {\it poles} (at which the function
$\varphi$ does not vanish) and {\it ambiguous points} (at which
$\varphi=0$).

One can also define a divisor of a meromorphic function. Denote by
$A=N\cup P$ the union of the set $N$ of zeros and polar locus $P$
of the function $f$. Any regular point $a^0$ of the set $A$ can
belong to only one irreducible component of $N$ or $P$. Thus, one
can define the {\it order} of this component as the degree (taken
with "+", if $a^0\in N$, and with "$-$", if $a^0\in P$) of the
corresponding factor in the decomposition of the function
$\varphi$ or $\psi$ into irreducible factors. Then the {\it
divisor} of the meromorphic function $f$ is the pair $(A,
\kappa)$, where $\kappa=\kappa(a)$ is an integer-valued function
on the set of regular points of $A$ (which takes a constant value
on each its irreducible component, this value is equal to the
order of a component). The pair $(P,\kappa)$ is called the {\it
polar divisor} of the meromorphic function $f$. By $(f)_{\infty}$
we will mean the restriction of $\kappa$ on regular points of $P$.

Let us return to the Schlesinger equation. The polar locus
$\Theta\subset Z$ of the extended matrix functions
$B_1(a),\ldots,B_n(a)$ is called the {\it Malgrange
$\Theta$-divisor}\footnote{In view of the above definition of a
divisor, here the term "divisor" is not precise enough.}. If we
consider the system (\ref{syst0}) as an equation for horizontal
sections of the logarithmic connection $\nabla_0$ (with
singularities $a_1^0,\ldots,a_n^0$) in the holomorphically trivial
vector bundle $E_0$ of rank $p$, then the set $\Theta$ corresponds
to those points $a^*$, where the bundle $E_{a^*}$ associated to
the parameter $a^*$ in the isomonodromic deformation $(E_a,
\nabla_a)_{a\in Z}$ of $(E_0, \nabla_0)$ is not holomorphically
trivial (see details in \cite{Ma}).

In what follows we will use the theorem of Bolibrukh \cite{Bo5},
\cite{Bo6} (the proof also can be found in \cite{Go}) describing a
general solution of the Schlesinger equation near the
$\Theta$-divisor in the case $p=2$. For the polar locus $P\subset
Z$ of a function $f$ meromorphic on $Z$, and $a^*\in P$, let us
denote by $\Sigma_{a^*}(f)$ the sum of orders of all irreducible
components of $P\cap D(a^*)$.
\medskip

{\bf Theorem 1.} {\it Let the monodromy of the two-dimensional
$(p=2)$ Schlesinger isomonodromic family $(\ref{fam})$ be
irreducible, $a^*$ an arbitrary point of $\Theta$ and
$E_{a^*}\cong{\cal O}(k)\oplus{\cal O}(-k)$. Then
$\Sigma_{a^*}(B_i)\geqslant-2k$, $i=1,\ldots,n$.}
\medskip

{\bf Remark 1.} As known \cite{Bo6}, $2k\leqslant n-2$. Thus, the
estimate of Theorem 1 can be written in the form
$\Sigma_{a^*}(B_i)\geqslant2-n$, furthermore
$\Sigma_{a^*}(B_i)\geqslant3-n$ in the case of odd $n$.
\medskip

The following auxiliary lemma will be used later.
\medskip

{\bf Lemma 1.} {\it Consider a two-dimensional Schlesinger
isomonodromic family of the form
$$
\frac{dy}{dz}=\left(\sum_{i=1}^n\frac{B_i(a)}{z-a_i}\right)y,\qquad
\sum_{i=1}^nB_i(a)=K={\rm diag}(\theta,-\theta),\quad
\theta\in{\mathbb C},
$$
and the function $b(a)=\sum_{i=1}^nb_i^{12}(a)a_i$, where
$b_i^{12}(a)$ are the upper-right elements of the matrices
$B_i(a)$ respectively. Then the differential of the function
$b(a)$ is given by the formula}
$$
db(a)=(2\theta+1)\sum_{i=1}^nb_i^{12}(a)da_i.
$$

{\bf Proof.} The differential $db(a)$ has the form
$$
db(a)=\sum_{i=1}^na_idb_i^{12}(a)+\sum_{i=1}^nb_i^{12}(a)da_i.
$$
To find the first of the two latter summands, let us use the
Schlesinger equation for the matrices $B_i(a)$. Then we have
\begin{eqnarray*}
\sum_{i=1}^n a_i\,dB_i(a)&=&-\sum_{i=1}^n\sum_{j=1, j\ne i}^n
a_i\frac{[B_i(a),B_j(a)]}{a_i-a_j}\,d(a_i-a_j)=
-\sum_{i=1}^n\sum_{j>i}^n\,[B_i(a),B_j(a)]d(a_i-a_j)=\\
&=&-\sum_{i=1}^n\Bigl[B_i(a),\sum_{j=1, j\ne i}^n B_j(a)\Bigr]da_i
=-\sum_{i=1}^n[B_i(a),K]da_i.
\end{eqnarray*}

The upper-right element of the latter matrix 1-form is equal to
$\sum_{i=1}^n2\theta b_i^{12}(a)da_i$, hence
$\sum_{i=1}^na_idb_i^{12}(a)=2\theta\sum_{i=1}^n b_i^{12}(a)da_i$,
and $db(a)=(2\theta+1)\sum_{i=1}^nb_i^{12}(a)da_i$. \hfill{$\Box$}
\\
\\

{\bf \S3. Schlesinger isomonodromic deformations and Garnier
systems}
\medskip

Let us recall the relationship between Schlesinger isomonodromic
deformations and Garnier systems.

Consider a two-dimensional Schlesinger isomonodromic family
\begin{eqnarray}\label{schl2n}
\frac{dy}{dz}=\left(\sum_{i=1}^{n+2}\frac{B_i(a)}{z-a_i}\right)y,
\qquad B_i(a^0)=B_i^0\in sl(2,\mathbb C),
\end{eqnarray}
of Fuchsian systems with singular points $a_1,\ldots,a_n$,
$a_{n+1}=0$, $a_{n+2}=1$, $a_{n+3}=\infty$ which depends
holomorphically on the parameter $a=(a_1,\ldots,a_n)\in D(a^0)$,
where $D(a^0)$ is a disk of small radius centered at the point
$a^0$ of the space $({\mathbb C}\setminus\{0,1\})^n
\setminus\bigcup_{i\ne j}\{a_i=a_j\}$. Denote by $\pm\beta_i$ the
eigenvalues of the matrices $B_i(a)$ respectively. Recall that the
isomonodromic deformation preserves the eigenvalues of the residue
matrices $B_i(a)$. As follows from the Schlesinger equation, the
matrix residue at the infinity is constant. We assume that it is a
diagonalisable matrix, i.~e.,
$\sum_{i=1}^{n+2}B_i(a)=-B_{\infty}={\rm diag}(-\beta_{\infty},
\beta_{\infty})$.

By Malgrange's theorem the matrix functions
$$
B_i(a)=\left(\begin{array}{cc} b_i^{11}(a) & b_i(a) \\
                               b_i^{21}(a) & b_i^{22}(a)
              \end{array}\right)
$$
can be extended to the universal cover $Z'$ of the space
$({\mathbb C}\setminus\{0,1\})^n\setminus\bigcup_{i\ne
j}\{a_i=a_j\}$ as meromorphic functions (holomorphic off the
analytic subset $\Theta$ of codimension one).

Denote by $B(z,a)$ the coefficient matrix of the family
(\ref{schl2n}). Since the upper-right element of the matrix
$B_{\infty}$ equals zero, for every fixed $a$ the same element of
the matrix $z(z-1)(z-a_1)\ldots(z-a_n)B(z,a)$ is a polynomial
$P_n(z,a)$ of degree $n$ in $z$. We denote by
$u_1(a),\ldots,u_n(a)$ the roots of this polynomial and define the
functions $v_1(a),\ldots,v_n(a)$:
$$
v_j(a)=\sum_{i=1}^{n+2}\frac{b_i^{11}(a)+\beta_i}{u_j(a)-a_i},
\qquad j=1,\ldots,n.
$$
Then the following statement takes place: {\it the pair
$(u(a),v(a)) =(u_1,\ldots,u_n,v_1,\ldots,v_n)$ satisfies the
Garnier system $(\ref{garnier})$ with the parameters
$2\beta_1,\ldots,2\beta_{n+2},2\beta_{\infty}-1$} (see proof of
Proposition 3.1 from \cite{Ok}, or \cite{GP}, Cor. 6.2.2 (p.
207)).

One can express the coefficients of the polynomial $P_n(z,a)$ in
terms of the upper-right elements $b_i(a)$ of the matrices
$B_i(a)$. Let
$$
\sigma_1(a)=\sum_{i=1}^{n+2}a_i,\quad \sigma_2(a)=\sum_{1\leqslant
i<j\leqslant n+2}a_ia_j,\quad\ldots,\quad
\sigma_{n+1}(a)=a_1\ldots a_n
$$
be the elementary symmetric polynomials in $a_1,\ldots,a_n$,
$a_{n+1}=0$, $a_{n+2}=1$, and $Q(z)=\prod_{i=1}^{n+2}(z-a_i)$.
Then
$$
P_n(z,a)=\sum_{i=1}^{n+2}b_i(a)\frac{Q(z)}{z-a_i}=:
b(a)z^n+f_1(a)z^{n-1}+\ldots+f_n(a)
$$
(recall that $\sum_{i=1}^{n+2}b_i(a)=0$). By the Vi\`ete theorem
one has
\begin{eqnarray*}
b(a)&=&\sum_{i=1}^{n+2}b_i(a)(-\sigma_1(a)+a_i)=
\sum_{i=1}^{n+2}b_i(a)a_i=\sum_{i=1}^nb_i(a)a_i+
b_{n+2}(a),\\
f_1(a)&=&\sum_{i=1}^{n+2}b_i(a)\Bigl(\sigma_2(a)-\sum_{j=1,j\ne
i}^{n+2}a_ia_j\Bigr)=-\sum_{1\leqslant i<j\leqslant n+2}
(b_i(a)+b_j(a))a_ia_j.
\end{eqnarray*}

In the similar way,
$$
f_k(a)=(-1)^k\sum_{1\leqslant i_1<\ldots<i_{k+1}\leqslant n+2}
(b_{i_1}(a)+\ldots+b_{i_{k+1}}(a))a_{i_1}\ldots a_{i_{k+1}}
$$
for each $k=1,\ldots,n$.

Alongside formulae for the transition from a two-dimensional
Schlesinger isomonodromic family with $sl(2,{\mathbb C})$-residues
to a Garnier system, there also exist formulae for the inverse
transition (see \cite{Ok}, Prop. 3.2). This allows to suggest some
addition to Garnier's theorem (which claims that the elementary
symmetric polynomials $F_i(a)=\sigma_i(u_1(a),\ldots,u_n(a))$ of
solutions of a Garnier system are meromorphic on $Z'$).

By the {\it linear monodromy} of a solution of a Garnier system we
will mean the monodromy of the corresponding two-dimensional
Schlesinger isomonodromic family.
\medskip

{\bf Theorem 2.} {\it Let $(u(a),v(a))$ be a solution of the
Garnier system $(\ref{garnier})$ that has an irreducible linear
monodromy, and $\Delta_i$ denotes the polar locus of the function
$F_i$, $i=1,\ldots,n$. Then

{\rm a)} in the case $\theta_{\infty}=0$ and $u_i(a)\not\equiv
u_j(a)$ for $i\ne j$, one has $\Sigma_{a^*}(F_i)\geqslant-n-1$ for
any point $a^*\in\Delta_i$;

{\rm b)} in the case $\theta_{\infty}\ne0$ one has
$\Sigma_{a^*}(F_i)\geqslant-n$ for any point $a^*\in\Delta_i$, may
be, with the exception of some subset $\Delta^0\subset\Delta_i$ of
positive codimension $($in any case $(F_i)_{\infty}\geqslant-n)$.}
\medskip

{\bf Proof.} Consider the family (\ref{schl2n}) with the
irreducible monodromy corresponding to the given solution, and the
functions $b(a)$, $f_1(a),\ldots,f_n(a)$ constructed by the
residue matrices $B_i(a)$. By the Vi\`ete theorem,
$F_i(a)=(-1)^if_i(a)/b(a)$. Due to Theorem 1 and Remark 1, for
each function $f_i$ and any point $a^*$ of the $\Theta$-divisor of
the family (\ref{schl2n}) one has
$\Sigma_{a^*}(f_i)\geqslant-n-1$.

By Lemma 1 we have $db(a)=-\theta_{\infty}\sum_{i=1}^nb_i(a)da_i$,
where $\theta_{\infty}=2\beta_{\infty}-1$.

a) In the case $\theta_{\infty}=0$ one has $db(a)\equiv0$ for all
$a\in Z'$, hence $b(a)\equiv{\rm const}\ne0$. Indeed, if
$b(a)\equiv0$, then $P_n(z,a)$ is a polynomial of degree $n-1$ in
$z$, and $u_i(a)\equiv u_j(a)$ for some $i\ne j$, which
contradicts the conditions of the theorem. Thus,
$\Sigma_{a^*}(F_i)=\Sigma_{a^*}(f_i)\geqslant-n-1$ in this case.

b) In the case $\theta_{\infty}\ne0$
\begin{eqnarray}\label{formulae}
b_i(a)&=&-\frac1{\theta_{\infty}}\frac{\partial b(a)}{\partial
a_i},\qquad
i=1,\ldots,n; \nonumber\\
b_{n+2}(a)&=&b(a)-\sum_{i=1}^nb_i(a)a_i, \quad
b_{n+1}(a)=-b_{n+2}(a)-\sum_{i=1}^nb_i(a).
\end{eqnarray}
Thus, if the function $b$ is holomorphic at a point $a'\in Z'$, so
are the functions $b_i$, $i=1,\ldots,n+2$, and hence, the
functions $f_i$. Therefore, the points $a^*\in\Delta_i$ can be of
two types: such that $b(a^*)=0$ (then
$\Sigma_{a^*}(F_i)\geqslant-1$, since the function $b$ is
irreducible\footnote{Indeed, if for some $a'\in\{b(a)=0\}$ one has
$db(a')\equiv0$, then $\sum_{i=1}^nb_i(a')da_i\equiv0$ and
$b_1(a')=\ldots=b_n(a')=0$. Taking into consideration the
relations (\ref{formulae}), one gets also $b_{n+2}(a')=0$ and
$b_{n+1}(a')=0$. This contradicts the irreducibility of the
monodromy of the family (\ref{schl2n}).}) or that belong to the
polar locus $\Delta\subset\Theta$ of the function $b$.

Denote by $\Delta^0\subset\Delta$ the set of ambiguous points of
the function $b$. Then in a neighbourhood of any point
$a^*\in\Delta\setminus\Delta^0$ it can be presented in the form
\begin{eqnarray}\label{divisor1}
b(a)=\frac{h(a)}{\tau_1^{j_1}(a)\ldots\tau_r^{j_r}(a)}, \quad
j_1\geqslant1,\ldots,j_r\geqslant1,
\end{eqnarray}
where the functions $\tau_i$, $h$ are holomorphic near $a^*$,
$h(a^*)\ne0$, furthermore $\tau_i$ are irreducible at $a^*$, just
as
\begin{eqnarray}\label{divisor2}
f_i(a)=\frac{g(a)}{\tau_1^{k_1}(a)\ldots\tau_r^{k_r}(a)}, \quad
k_1+\ldots+k_r\leqslant n+1,
\end{eqnarray}
where the function $g$ is holomorphic near $a^*$. Thus,
$$
\frac{f_i(a)}{b(a)}=
\frac{g(a)}{\tau_1^{k_1}(a)\ldots\tau_r^{k_r}(a)}:
\frac{h(a)}{\tau_1^{j_1}(a)\ldots\tau_r^{j_r}(a)}=
\frac{g(a)/h(a)}{\tau_1^{k_1-j_1}(a)\ldots\tau_r^{k_r-j_r}(a)},
$$
therefore,
$$
\Sigma_{a^*}(F_i)=-\sum_{\alpha}(k_{\alpha}-j_{\alpha})\geqslant-n
$$
(the sum is taken with respect to such indices $\alpha$ that
$k_{\alpha}-j_{\alpha}>0$), which proves the first part of the
statement b).

In a neighbourhood of a point $a^*\in\Delta^0$ the decompositions
(\ref{divisor1}), (\ref{divisor2}) take place for the functions
$b$, $f_i$  respectively, but $h(a^*)=0$. However, due to the
irreducibility of $b$, all irreducible factors of $h$ in its
decomposition $h(a)=h_1(a)\ldots h_s(a)$ near $a^*$ are distinct
(we can assume also that none of $h_i$ coincides with some of
$\tau_l$). Since $k_l-j_l\leqslant n$ for all $l=1,\ldots,r$, the
second part of the statement b) follows from the decomposition
$$
\frac{f_i(a)}{b(a)}=\frac{g(a)}{h_1(a)\ldots h_s(a)\,
\tau_1^{k_1-j_1}(a)\ldots\tau_r^{k_r-j_r}(a)}.
$$
{\hfill $\Box$}
\medskip

{\bf Remark 2.} As follows from Remark 1, in all estimates of
Theorem 2 one can substitute $n$ by $n-1$ in the case of even $n$.

In particular, the polar loci of the functions
$F_1(a)=u_1(a)+u_2(a)$ and $F_2(a)=u_1(a)u_2(a)$, where
$(u_1,u_2,v_1,v_2)$ is a solution of the Garnier system ${\cal
G}_2(\theta_1,\ldots,\theta_4, \theta_{\infty})$ corresponding to
a two-dimensional Schlesinger isomonodromic family with five
singular points and irreducible monodromy, are analytical
submanifolds with $(F_i)_{\infty} \geqslant-2$. (Note that a
bundle $E_{a^*}$ corresponding to a point $a^*$ of the
$\Theta$-divisor of this family has the form $E_{a^*}\cong{\cal
O}(1)\oplus{\cal O}(-1)$, which implies the regularity of the
$\Theta$-divisor; see \cite{Bo3}, \cite{Go}.)
\medskip

M.\,Mazzocco \cite{Maz2} has shown that the solutions of the
Garnier system (\ref{garnier}), that have {\it reducible} linear
monodromy, are classical functions (in each variable, in sense of
Umemura \cite{Um}) and can be expressed via solutions of
Lauricella hypergeometric equations. We discuss this case in more
details in the following section.
\\
\\

{\bf \S4. Garnier systems and Lauricella hypergeometric equations}
\medskip

Consider the {\it Lauricella hypergeometric equation}
$E_D(\alpha,\beta_1,\ldots,\beta_n,\gamma)$
\begin{eqnarray*}
a_i(1-a_i)\frac{\partial^2u}{\partial a_i^2}&+&(1-a_i)
\sum_{j=1}^na_j\frac{\partial^2u}{\partial a_i
\partial a_j}+(\gamma-(\alpha+\beta_i+1)a_i)\frac{\partial u}{\partial
a_i}-\\&-& \beta_i\sum_{j=1}^na_j\frac{\partial u}{\partial
a_j}-\alpha\beta_iu=0, \qquad i=1,\ldots,n; \\
(a_i-a_j)\frac{\partial^2u}{\partial a_i
\partial a_j}&+&\beta_i\frac{\partial u}{\partial
a_j}-\beta_j\frac{\partial u}{\partial a_i}=0,\qquad
i,j=1,\ldots,n.
\end{eqnarray*}

This is a system of linear partial differential equations of
second order for a complex function $u(a_1,\ldots,a_n)$, where
$\alpha,\beta_1,\ldots,\beta_n,\gamma\in{\mathbb C}$ are
constants. The system is defined on the space $B=({\mathbb
C}\setminus\{0,1\})^n\setminus\bigcup_{i\ne j}\{a_i=a_j\}$. As
shown in \cite{GP} (proof of Proposition 9.1.4, p.249), the
vector-function $y(a)=(u,a_1\frac{\partial u}{\partial
a_1},\ldots, a_n\frac{\partial u}{\partial a_n})^{\top}$,
$a=(a_1,\ldots,a_n) \in B$, satisfies a completely integrable
linear Pfaffian system
\begin{eqnarray}\label{pfaff}
dy=\omega y,\qquad y(a)\in{\mathbb C}^{n+1}.
\end{eqnarray}
Therefore, the set of solutions of a Lauricella hypergeometric
equation near every point $a\in B$ forms an $(n+1)$-dimensional
vector space, and solutions can be extended holomorphically to the
universal cover $Z'$ of the space $B$.

Further we recall how particular solutions of certain Garnier
systems are connected with Lauricella hypergeometric equations and
study movable singularities of such solutions.

As we noted, for $n>1$ the Garnier system (\ref{garnier})
generically does not satisfy the Painlev\'e property (coordinates
$u_1,\ldots,u_n$ of its solution are defined as roots of a
polynomial of degree $n$), but it can be transformed by a certain
(symplectic) transformation $(a,u,v,H)\mapsto(s,q,p,\widetilde
H)$, $\sum_{i=1}^n(p_idq_i-\widetilde H_i ds_i)
=\sum_{i=1}^n(v_idu_i-H_ida_i)$, into a Hamiltonian system ${\cal
H}_n(\theta)$
$$
\frac{\partial q_j}{\partial s_i}=\frac{\partial\widetilde
H_i}{\partial p_j}, \qquad \frac{\partial p_j}{\partial s_i}=
-\frac{\partial\widetilde H_i}{\partial q_j}, \qquad
i,j=1,\ldots,n,
$$
satisfying the Painlev\'e property (see \cite{GP}, Ch. III, \S7).
Recall this transformation introducing the functions
\begin{eqnarray}\label{funct}
M_i(a,u)=-\frac{(a_i-u_1)\ldots(a_i-u_n)}{\prod_{j=1,j\ne
i}^{n+2}(a_i-a_j)}, &  &
M^{k,i}(a,u)=\frac{u_k(u_k-1)(u_k-a_1)\ldots(u_k-a_n)}{(u_k-a_i)\prod_{j=1,j\ne
k}^n(u_k-u_j)}, \nonumber \\
 &  & i, k=1,\ldots,n.
\end{eqnarray}

Then the transformation $(a,u,v)\mapsto(s,q,p)$ is given by the
formulae
\begin{eqnarray}\label{transf}
s_i=\frac{a_i}{a_i-1}, \quad q_i=-a_iM_i, \quad
p_i=(1-a_i)\sum_{k=1}^n\frac{M^{k,i}\,v_k}{u_k(u_k-1)}, \qquad
i=1,\ldots,n,
\end{eqnarray}
furthermore
\begin{eqnarray}\label{transf2}
v_i=\sum_{k=1}^n\frac{q_kp_k}{u_i-a_k}, \qquad i=1,\ldots,n,
\end{eqnarray}
while the new Hamiltonians
$$
\widetilde H_i=-(1-a_i)^2\Bigl(H_i+\sum_{j=1}^np_j\frac{\partial
q_j}{\partial a_i}\Bigr)=
\frac1{s_i(s_i-1)}\Bigl(\sum_{j,k=1}^nE^i_{jk}(s,q)p_jp_k-
\sum_{j=1}^nF^i_j(s,q)p_j+\varkappa q_i\Bigr)
$$
are polynomial in $(q, p)=(q_1,\ldots,q_n,p_1,\ldots,p_n)$. Here
$E^i_{jk}(s,q)$, $F^i_j(s,q)$ are polynomials in $q$ of third and
second degrees respectively, and
$\varkappa=\frac14\left((\sum_{i=1}^{n+2}\theta_i-1)^2-
\theta_{\infty}^2\right)$.

In the case $\varkappa=0$ (i. e., $\sum_{i=1}^{n+2}\theta_i-1=
\pm\theta_{\infty}$) the system ${\cal H}_n(\theta)$ has solutions
of the form $(q,0)$, where $q$ is a solution of the system
\begin{eqnarray}\label{ric}
s_i(s_i-1)\frac{\partial q_j}{\partial s_i}=-F^i_j(s,q), \qquad
i,j=1,\ldots,n.
\end{eqnarray}
Since the right side of this system consists of polynomials of
second degree in $q$, the system may be considered as a
several-variables generalization of the classical Riccati
equation. Similarly to the classical case, the system (\ref{ric})
can be linearized by a suitable change of unknown. Exactly,
according to \cite{GP} (Th. 9.2.1, p.252), a solution
$(q_1,\ldots,q_n)$ of (\ref{ric}) can be presented in the form
$$
q_i(s)=\frac{s_i(s_i-1)}{\sum_{i=1}^{n+2}\theta_i-1}\left(\frac{\theta_i}
{s_i-1}+\frac1{f}\frac{\partial f}{\partial s_i}\right), \qquad
i=1,\ldots,n,
$$
where $f(s)$ is an arbitrary solution of the Lauricella
hypergeometric equation $E_D(1-\theta_{n+2}, \theta_1,\ldots,$
$\theta_n,\sum_{i=1}^{n+1}\theta_i)$.

The function $f$ is irreducible at its zeros (if $f(s^*)=0$ and
$\frac{\partial f}{\partial s_1}(s^*)=\ldots=\frac{\partial f}
{\partial s_n}(s^*)=0$, then $f\equiv0$ due to the uniqueness of
solution of the system (\ref{pfaff})), therefore a solution
$(q_1,\ldots,q_n)$ of the system (\ref{ric}), as well as any
linear combination $Q$ of $q_i$ (with holomorphic coefficients),
are meromorphic on $Z'$ and the polar locus of $Q$ is an
analytical submanifold with $(Q)_{\infty}=-1$.

Now consider a solution $(u,0)$ of the Garnier system
(\ref{garnier}) corresponding to a solution $(q,0)$ of the system
${\cal H}_n(\theta)$, with $\sum_{i=1}^{n+2}\theta_i-1=
\pm\theta_{\infty}$ (note that $v=0\Longleftrightarrow p=0$, by
(\ref{transf}) and (\ref{transf2})).  The elementary symmetric
polynomials $F_i(a)=\sigma_i(u_1(a),\ldots,u_n(a))$ are expressed
via linear combinations of $q_i$ with holomorphic coefficients.
Indeed, let
$$
Q_i(a):=\prod_{j=1,j\ne i}^{n+2}(a_i-a_j),\qquad i=1,\ldots,n+2.
$$
Then, as follows from the formulae (\ref{funct}) and
(\ref{transf}),
$$
a_i^n-F_1(a)a_i^{n-1}+\ldots+(-1)^nF_n(a)=\frac{Q_i(a)}{a_i}q_i,
\qquad i=1,\ldots,n.
$$
Thus, the vector $(-F_1,\ldots,(-1)^nF_n)$ is a solution of the
system of linear equations with the coefficients matrix whose
determinant is the Vandermonde determinant.

The above reasonings lead to the following statement.
\medskip

{\bf Proposition 1.} {\it Let $(u(a), 0)$ be a solution of the
Garnier system $(\ref{garnier})$ with $\sum_{i=1}^{n+2}\theta_i-1=
\pm\theta_{\infty}$ $($this solution has the reducible linear
monodromy$)$. Then the polar loci of the functions $F_i$ are
analytical submanifolds and $(F_i)_{\infty}=-1$}.
\medskip

M.\,Mazzocco \cite{Maz2} has shown that any solution of the
Garnier system  (\ref{garnier}) having {\it reducible} linear
monodromy can be expressed classically via particular solutions
from the above proposition. We hope that a careful reading of her
article will lead to the obtaining of concluding estimates for the
elementary symmetric polynomials depending on an {\it arbitrary}
solution of (\ref{garnier}).

\bigskip

Institute for Information Transmission Problems of the Russian
Academy of Sciences,

Bolshoy Karetny per. 19, Moscow, 127994, Russia
\\

rgontsov@inbox.ru

\end{document}